\documentclass[11pt]{amsart}

\usepackage[utf8]{inputenc}
\usepackage[english]{babel}
\usepackage{amsmath,amsfonts,amssymb,latexsym,amsthm}
\usepackage{hyperref}
\hypersetup{colorlinks} 

\hypersetup{
colorlinks=true,
linkcolor=blue,
filecolor=magenta,
urlcolor=cyan,
}

\usepackage{caption}
\usepackage[mode=buildnew]{standalone}

\usepackage{amssymb,verbatim}
\usepackage{amsmath,amsfonts}
\usepackage[mathscr]{euscript} 
\usepackage{amsthm}
\usepackage{url}
\usepackage{graphicx} 
\usepackage{enumitem} 
\usepackage{bbm} 
\usepackage{bm} 
\usepackage{mathrsfs} 
\usepackage{cancel} 

\usepackage[left=1in, right=1in, top=1in, bottom=1in]{geometry}
\usepackage{graphicx}
\usepackage{tikz}

\renewcommand{\P}{\mathbb{P}}
\newcommand{\Pbb}{\mathbb{P}^{\textrm{bb}}}
\newcommand{\Pbbr}{\mathbb{P}^{\text{\em bb}}}

\newcommand{\PWZ}{\mathbb{P}^{\textrm{wz}}}
\newcommand{\PWZr}{\mathbb{P}^{\text{\em wz}}}

\newcommand{\Palt}{\mathbb{P}^{\textrm{alt}}}
\newcommand{\Paltr}{\mathbb{P}^{\text{\em alt}}}

\usepackage[colorinlistoftodos]{todonotes}

\RequirePackage{cleveref}
\usepackage{hypcap}
\hypersetup{colorlinks=true, citecolor=darkblue, linkcolor=darkblue}
\definecolor{darkblue}{rgb}{0.0,0,0.7}
\newcommand{\darkblue}{\color{darkblue}}

\definecolor{darkred}{rgb}{0.68,0,0}
\newcommand{\darkred}{\color{darkred}}

\definecolor{darkgreen}{rgb}{0,.38,0}

\definecolor{magenta}{rgb}{.51, 0, .51}

\newcommand{\defn}[1]{\emph{\darkblue #1}}
\newcommand{\defna}[1]{\emph{\darkred #1}}

\setlist[enumerate]{
	label=\textnormal{({\roman*})},
	ref={\roman*}}

\makeatletter
\def\th@plain{%
	\thm@notefont{}
	\itshape 
}
\def\th@definition{%
	\thm@notefont{}
	\normalfont 
}
\makeatother

\newtheorem{thm}{Theorem}[section]
\newtheorem{lemma}[thm]{Lemma}

\newtheorem*{claim*}{Claim}

\newtheorem{prop}[thm]{Proposition}
\newtheorem{conj}[thm]{Conjecture}
\newtheorem{conv}[thm]{Convention}

\theoremstyle{definition}

\newtheorem{rem}[thm]{Remark}

\numberwithin{figure}{section}
\numberwithin{equation}{section}

\def\zz{\mathbb Z}

\def\ov{\overline}

\def\sm{\smallsetminus}

\def\Om{\Omega}

\def\ga{\gamma}

\def\vk{\varkappa}

\def\cC{\mathcal C}

\def\cA{\mathcal A}

\def\ssu{\subset}

\def\<{\langle}
\def\>{\rangle}

\def\Ups{{\small {\Upsilon}}}

\def\0{{\mathbf 0}}

\def\.{\hskip.06cm}
\def\ts{\hskip.03cm}

\def\lra{\leftrightarrow}

\def\La{\Lambda}

\def\.{\hskip.06cm}
\def\ts{\hskip.03cm}

\def\nin{\noindent}

\def\sz{{}}
\def\so{{\prime}}

\newcommand{\textsu}[1]{\textup{\textsf{#1}}}

\newcommand{\ComCla}[1]{\textup{\textsu{#1}}}

\newcommand{\sharpP}{\ComCla{\#P}}
\newcommand{\SP}{\ComCla{\#P}}

\def\SP{\sharpP}

\title{The bunkbed conjecture is false}
\date{\today}

 \author{Nikita Gladkov}
 \address[Nikita Gladkov]{Department of Mathematics, UCLA,  Los Angeles, CA 90095.}
 \email{\texttt{gladkovna@math.ucla.edu}}

 \author[\ts Igor Pak]{Igor Pak}
 \address[Igor Pak]{Department of Mathematics, UCLA,  Los Angeles, CA 90095.}
 \email{\texttt{pak@math.ucla.edu}}

 \author{Aleksandr Zimin}
 \address[Aleksandr Zimin]{Department of Mathematics, MIT, Cambridge, MA 02139.}
 \email{\texttt{azimin@mit.edu}}

\begin{document}

\begin{abstract}
We give an explicit counterexample to the \emph{bunkbed conjecture}
introduced by Kasteleyn in 1985.  The counterexample is given by a
planar graph on $7222$ vertices, and is built on the recent work
of Hollom (2024).
\end{abstract}
\maketitle

\section{Introduction}

The \emph{bunkbed conjecture} (BBC) is a celebrated open problem in probability
introduced by Kasteleyn in 1985, see  \cite[Remark 5]{BKa}.
The conjecture is both natural and intuitively obvious, but has
defied repeated proof attempts; it is known only in a few special cases.
In this paper we disprove the conjecture without resorting to
computer experiments (cf.\ Section~\ref{s:experimental}).

\smallskip

Let $G=(V,E)$ be a connected graph, possibly infinite and with multiple edges.
In \emph{Bernoulli bond percolation}, each edge is \emph{deleted}
independently at random with probability $1-p$,
and otherwise \emph{retained} with probability \ts $p\in [0,1]$.  Equivalently,
this model gives a random subgraph of $G$ weighted by the number of edges.
For $p=\frac12$ we obtain a uniform random subgraph of~$G$.
See \cite{BR06,Gri99} for standard results and \cite{Dum18,Wer09}
for recent overview of percolation.

Let \ts $\P_p[u \lra v]$ \ts denote the probability that vertices
$u,v\in V$ are connected.  It is often of interest to compare these
probabilities, as computing them exactly is $\SP$-hard \cite{PB83}.
For example, the classical \emph{Harris--Kleitman inequality}, a special
case of the \emph{FKG inequality}, implies that \ts
$\P_p[u \lra v]\le \P_p[u \lra v \. | \. u \lra w]$ \ts for all \ts $u,v,w\in V$,
see e.g. \cite[Ch.~6]{AS16}.
Harris used this to prove that  the
\emph{critical probability} \ts $p_c(G) :=\inf\{p \. : \. \P_p(G)>0\}$ \ts
satisfies \ts $p_c(\zz^2)\ge \frac12$ \ts \cite{Har60},
in the first step towards Kesten's remarkable exact value
\ts $p_c(\zz^2)= \frac12$ \ts \cite{Kes80}.  Considerations of
percolation monotonicity on \ts $\zz^2$ (see~$\S$\ref{ss:finrem-2PF}),
led Kasteleyn to the following problem.

\smallskip

Fix a finite connected graph \ts $G=(V,E)$ \ts and a subset \ts $T\subseteq V$.
A \emph{bunkbed graph} \ts $\ov G=(\ov V, \ov E)$ \ts is a subgraph of the graph product
\ts $G \times K_2$ \ts defined as follows.  Take two copies of $G$, which
we denote $G$ and $G'=(V',E')$, and add all edges of the form $(w,w')$,
where $w\in T$ and $w'$ is a corresponding vertex in~$T'$;
we denote this set of edges by~$\ov T$.  The resulting bunkbed graph
has \ts $\ov V = V\cup V'$ \ts and \ts $\ov E = E \cup E' \cup \ov T$.

In the \defn{bunkbed percolation}, the usual bond percolation is performed only on edges
in $G$ and $G'$, while all edges in $\ov T$ are retained (i.e., not deleted).  We use
\ts $\Pbb_p[u \lra v]$ \ts to denote connecting probability in this case.
The vertices in $T$ are called \emph{transversal} and the edges in $\ov T$ are called
\emph{posts}, to indicate their special status.  See e.g.\ \cite{Lin11, RS},
for these and several other equivalent models of the bunkbed percolation.
We refer also to \cite[$\S$4.1]{Gri23}, \cite[$\S$5.5]{Pak-OPAC} and \cite{Rud21}
for recent overviews and connections to other areas.

\smallskip

\begin{conj}[{\rm \defna{bunkbed conjecture}}{}] \label{conj:BB}
Let \ts $G=(V,E)$ \ts be a connected graph, let \ts $T\subseteq V$, and
let \ts $0< p <1$.  Then, for all \ts $u,v\in V$, we have:
$$\Pbbr_p[\ts u \lra v \ts ] \. \ge \. \Pbbr_p[\ts u \lra v' \ts ]\ts.
$$
\end{conj}

\smallskip

The bunkbed conjecture is known in a number of special cases, including wheels \cite{Lea09},
complete graphs \cite{dB16,dB18,HL19}, complete bipartite graphs \cite{Ric22},
and graphs symmetric w.r.t.\ the \ts $u \lra v$ \ts automorphism \cite{Ric22}.
It is also known for one \cite[Lemma~2.4]{Lin11} and for two transversal vertices \cite[$\S$6.3]{Lohr},
see also \cite{BKa,GZ+}.  Finally, the conjecture was recently proved in the \ts $p\uparrow 1$ \ts limit \cite{HNK21,Hol24}.
\smallskip

\begin{thm}\label{thm:main}
There is a connected planar graph \ts $G=(V,E)$ \ts with  \ts $|V|= 7222$ \ts vertices and \ts
$|E|=14442$ \ts edges,
a subset \ts $T\ssu V$ \ts with three transversal vertices, and vertices \ts $u,v \in V$, s.t.
\[\Pbbr_{\frac12}[\ts u \lra v \ts ] \. < \. \Pbbr_{\frac12}[ \ts u \lra v' \ts ]\ts.
\]
In particular, the bunkbed conjecture is false.
\end{thm}

\smallskip
The result is surprising since analogous inequalities for simple random walks
and for the Ising model on bunkbed graphs were proved by H\"{a}ggstr\"{o}m \cite{Hag98,Hag03},
cf.~$\S$\ref{ss:finrem-RC}.
Recall that three is the smallest number of transversal vertices we can have to disprove
the conjecture.  On the other hand, the total number of vertices is unlikely to be optimal,
see Remark~\ref{r:counter-example} and Section~\ref{s:experimental}.

The proof of the theorem is based on an example of Hollom \cite{H} refuting
the $3$-uniform hypergraph version of the BBC.  Unfortunately, Hollom's example
alone cannot disprove the conjecture since it is impossible
to find a gadget graph simulating a single $3$-hyperedge using bond percolation
 \cite[Thm~1.5]{GZ}.

We give a robust version of Hollom's construction using the approach
in \cite{Gla24,GZ}.  The proof of Theorem~\ref{thm:main} occupies most
of the paper.  It is self-contained modulo Hollom's result which is small
enough to be checked by hand.  In Section~\ref{s:complete}, we extend the
theorem to the case when the set of transversal vertices is not fixed but
chosen uniformly at random from~$V$, see Theorem~\ref{thm:complete}.
We conclude with discussion of our computer experiments in
Section~\ref{s:experimental}, and final remarks in Section~\ref{s:finrem}.

\medskip

\section{Notation} \label{s:notation}

In percolation, deleted edges are called \defn{closed} \ts while retained edges
are called \defn{open}.  Note that there are several different models of percolation
and variations on the bunkbed conjecture (BBC), see $\S$\ref{ss:finrem-models}.

A hypergraph is a collection of subsets of vertices;
to simplify the notation we use the same letter to denote both.
The hypergraph is called \defn{uniform} \ts if all hyperedges have
the same size. A \defn{path} in a hypergraph is a sequence \. $(v_0 \to v_1 \to \. \ldots \. \to v_\ell)$ \.
of vertices, such that \ts $v_{i-1}, v_i$ \ts lie in the same hyperedge, for
all \ts $1\le i \le \ell$.
We say that two vertices in a hypergraph are \defn{connected} \ts
if there is there is a path between them.
For further definitions and results on hypergraphs, see e.g.\
\cite[$\S$1.2]{Ber89}.

The notion of \defn{hypergraph percolation} \ts is
is a natural extension of graph percolation, and goes back to
the study of random hypergraphs, see e.g.\ \cite{SPS85}.
In recent years, the study of hypergraph percolation also
emerged in probabilistic and statistical physics literature,
see e.g.\ \cite{WZ} and \cite{BD24}, respectively.

\medskip

\section{Hypergraph percolation}\label{s:hyper}

\subsection{Hollom's example}\label{ss:hyper-ex}
%
%
%
Let $H$ be a finite connected hypergraph on the set $V$ of vertices.
We use \ts $\P_p[\ts u \lra v\ts]$ \ts to denote probability of
connectivity of vertices \ts $u,v\in V$ \ts in the
\defn{hypergraph percolation}, where each hyperedge $e$ in $H$
is retained with probability~$p$, or deleted with probability~$1-p$.

Let \ts $T\subseteq V$ \ts
be the set of transversal vertices.
Denote by \ts $\ov{H}$ \ts be the \defn{bunkbed hypergraph} with
levels \ts $H\simeq H'$, and \emph{vertical posts} which are the (usual) edges.
Note that \ts $\ov{H}$ \ts has \emph{horizontal} hyperedges and vertical posts.

In \cite{H}, Hollom considers the following natural hypergraph generalization
of the \emph{Alternative BBC}, see~$\S$\ref{ss:finrem-models}.
In the \defn{alternative bunkbed hypergraph percolation},
each hyperedge $e$ in $H$ is either deleted while the corresponding
hyperedge $e'$ in~$H'$ is retained with probability $\frac{1}{2}$, or
vice versa: hyperedge $h$ is retained and $h'$ is deleted.

\begin{figure}[hbt]
    \centering
    \begin{tikzpicture}[scale=4, thick]

    \definecolor{myyellow}{rgb}{0.98, 0.93, 0.36}  
    \colorlet{myyellow}{myyellow!40}
\definecolor{myblue}{rgb}{0.0, 0.0, 1.0}  
\definecolor{myred}{rgb}{0.8, 0.0, 0.0} 
\definecolor{mylightgray}{gray}{0.7}     

    \def\edgewidth{1mm}                     

    \node[circle, fill=mylightgray, draw, minimum size=8mm, label=center:$u_2$] (u2) at (0.5, 2) {};
    \node[circle, fill=mylightgray, draw, minimum size=8mm, label=center:$u_7$] (u7) at (1.5, 1.1) {};
    \node[circle, fill=mylightgray, draw, minimum size=8mm, label=center:$u_9$] (u9) at (2.5, 0.5) {};
    
    \node[circle, fill=white, draw, minimum size=8mm, label=center:$u_1$] (u1) at (1.5, 2.6) {};
    \node[circle, fill=white, draw, minimum size=8mm, label=center:$u_3$] (u3) at (2.5, 2) {};
    \node[circle, fill=white, draw, minimum size=8mm, label=center:$u_4$] (u4) at (0.5, 1.3) {};
    \node[circle, fill=white, draw, minimum size=8mm, label=center:$u_5$] (u5) at (1.25, 1.8) {};
    \node[circle, fill=white, draw, minimum size=8mm, label=center:$u_6$] (u6) at (1.75, 1.8) {};
    \node[circle, fill=white, draw, minimum size=8mm, label=center:$u_8$] (u8) at (0.5, 0.5) {};
    \node[circle, fill=white, draw, minimum size=8mm, label=center:$u_{10}$] (u10) at (1.5, 0) {};
    
    \fill[myyellow] (u1.center) -- (u2.center) -- (u3.center) -- cycle;
    \fill[myyellow] (u2.center) -- (u4.center) -- (u5.center) -- cycle;
    \fill[myyellow] (u3.center) -- (u6.center) -- (u9.center) -- cycle;
    \fill[myyellow] (u4.center) -- (u8.center) -- (u7.center) -- cycle;
    \fill[myyellow] (u7.center) -- (u5.center) -- (u6.center) -- cycle;
    \fill[myyellow] (u8.center) -- (u9.center) -- (u10.center) -- cycle;

    \draw[myblue, line width=\edgewidth] (u2) -- (u1);
    \draw[myblue, line width=\edgewidth] (u2) -- (u3);
    \draw[myblue, line width=\edgewidth] (u2) -- (u4);
    \draw[myblue, line width=\edgewidth] (u2) -- (u5);
    \draw[myblue, line width=\edgewidth] (u7) -- (u4);
    \draw[myblue, line width=\edgewidth] (u7) -- (u5);
    \draw[myblue, line width=\edgewidth] (u7) -- (u6);
    \draw[myblue, line width=\edgewidth] (u7) -- (u8);
    \draw[myblue, line width=\edgewidth] (u9) -- (u3);
    \draw[myblue, line width=\edgewidth] (u9) -- (u6);
    \draw[myblue, line width=\edgewidth] (u9) -- (u8);
    \draw[myblue, line width=\edgewidth] (u9) -- (u10);

    \draw[myred, line width=\edgewidth] (u1) -- (u3);
    \draw[myred, line width=\edgewidth] (u4) -- (u5);
    \draw[myred, line width=\edgewidth] (u4) -- (u8);
    \draw[myred, line width=\edgewidth] (u5) -- (u6);
    \draw[myred, line width=\edgewidth] (u3) -- (u6);
    \draw[myred, line width=\edgewidth] (u8) -- (u10);

    \node[circle, fill=mylightgray, draw, minimum size=8mm, label=center:$u_2$] (u2) at (0.5, 2) {};
    \node[circle, fill=mylightgray, draw, minimum size=8mm, label=center:$u_7$] (u7) at (1.5, 1.1) {};
    \node[circle, fill=mylightgray, draw, minimum size=8mm, label=center:$u_9$] (u9) at (2.5, 0.5) {};
    
    \node[circle, fill=white, draw, minimum size=8mm, label=center:$u_1$] (u1) at (1.5, 2.6) {};
    \node[circle, fill=white, draw, minimum size=8mm, label=center:$u_3$] (u3) at (2.5, 2) {};
    \node[circle, fill=white, draw, minimum size=8mm, label=center:$u_4$] (u4) at (0.5, 1.3) {};
    \node[circle, fill=white, draw, minimum size=8mm, label=center:$u_5$] (u5) at (1.25, 1.8) {};
    \node[circle, fill=white, draw, minimum size=8mm, label=center:$u_6$] (u6) at (1.75, 1.8) {};
    \node[circle, fill=white, draw, minimum size=8mm, label=center:$u_8$] (u8) at (0.5, 0.5) {};
    \node[circle, fill=white, draw, minimum size=8mm, label=center:$u_{10}$] (u10) at (1.5, 0) {};

\end{tikzpicture}
    \caption{Hollom's 3-uniform hypergraph $H$.}
    \label{fig:h}
\end{figure}
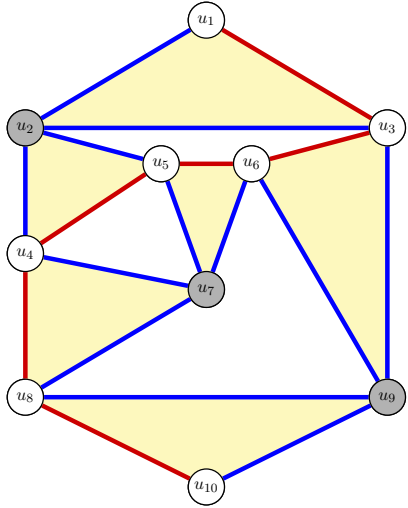

\begin{lemma}[{\rm Hollom \cite[Claim~5.1]{H}}{}]\label{lem:hollom}
Let $H$ be the hypergraph with six $3$-edges as in the Figure~\ref{fig:h},
and let \ts $T = \{u_2, u_7, u_9\}$ is the set of transversal vertices.
In the alternative bunkbed hypergraph percolation, we have:
$$
\Paltr[\ts u_1^{\sz}\lra u_{10}^{\so} \ts] \. = \. \frac{13}{64} \quad \text{and} \quad
\Paltr[\ts u_1^{\sz}\lra u_{10}^{\sz}\ts] \. = \. \frac{12}{64}\..
$$
\end{lemma}

We give a robust version of Hollom's construction.

\subsection{Robust hyperedge lemma}\label{ss:hyper-robust}
Note that in Hollom's example, each hyperedge has exactly one transversal vertex.
We explore this symmetry.

Consider the following \defn{WZ hypergraph percolation model} \ts introduced by
Wierman and Ziff in \cite{WZ} (see also \cite{GZ}).
Let \ts $e=\{a,b,c\}$ \ts be a hyperedge where
\ts $a$ \ts is a transversal vertex.
In the model, hyperedge $e$ is set to have

\smallskip

$\circ$ \ probability \ts $p_{abc}$ \ts to connect all three vertices,

$\circ$ \ probability \ts $p_{a|b|c}$ \ts to not connect any of the vertices,

$\circ$ \ probability \ts $p_{a|bc}$ \ts
to connect two non-transversal vertices, and

$\circ$ \ probability \ts $p_{ab|c}=p_{ac|b}$ \ts
to connect a transversal to a nontransversal vertex,

\smallskip

\nin
and these events are independent on all hyperedges.

\smallskip

Finally, we assume that these five probabilities sum up to~$1$:
$$p_{abc} \. + \. p_{a|b|c} \. + \. p_{a|bc} \. + \.
p_{ab|c} \. + \. p_{ac|b} \, = \. 1\ts.
$$
We say that vertices \ts $u$ \ts and \ts $v$ \ts are \defn{connected},
write \ts $u \lra v$, if they are connected in the hypergraph in a way
that every two vertices on a hyperedge are connected by the rules above.
We use \ts $\PWZ[\ts u \lra v\ts]$ \ts
to denote connection probabilities in this model, and drop the superscript
when the model is clear.

\begin{lemma}\label{lem:robust}
Let \ts $H$ \ts be Hollom's hypergraph in the Figure \ref{fig:h}, and let \ts
$T = \{u_2, u_7, u_9\}$ \ts be the set of transversal vertices.
Consider the WZ hypergraph percolation as described above, where the
connection probabilities satisfy
\begin{equation}\label{eq:390}
400 \. p_{a|bc} \, \le \, p_{abc} \. p_{a|b|c} \. - \. p_{ab|c}^2\ts.
\end{equation}
Then we have:
\begin{equation}\label{eq:semi}
\PWZr(u_1^{\sz}\lra u_{10}^{\sz}) \, < \, \PWZr(u_1^{\sz}\lra u_{10}^{\so})\ts.
\end{equation}
\end{lemma}

It was noted in \cite[Cor.~3.6]{Gla24}, that the RHS in \eqref{eq:390}
is always nonnegative:
\begin{equation}\label{eq:HK-abc}
p_{ab|c} \. p_{ac|b}  \. = \. p_{ab|c}^2  \. \le \. p_{abc} \. p_{a|b|c}\..
\end{equation}
This is a consequence of the Harris--Kleitman (HK)
inequality.  In fact, a slightly stronger inequality always holds (ibid.)
Since the LHS in \eqref{eq:390} is nonnegative, one can view this assumption
as strengthening the HK inequality in this case (cf.~$\S$\ref{ss:finrem-lemma}).

\smallskip

\subsection{Proof of Lemma~\ref{lem:robust}} \label{ss:hyper-proof}
Note that each configuration \ts $Q$ \ts of the
WZ~hypergraph percolation can be viewed
as a function \ts $\psi: H\cup H' \to \Ups$ \ts from hyperedges
in \ts $\ov{H}$ \ts to the set \ts $\Ups:= \{abc, ab|c, ac|b, a|bc, a|b|c\}$.
Here \ts $\psi$ \ts is the probability \ts $\P[Q]$ \ts of the configuration.
Clearly, for each hyperedge $e$ in~$H$, there are $25$ possibilities for \ts $\psi(e)$ \ts
and \ts $\psi(e')$.

To prove the lemma, we employ a combinatorial argument on all configurations that is similar
to the involution in the proofs of \cite[Lemma~2.3]{Lin11} and \cite[Lemma~2.5]{H}.
Denote by $\cC$ and $\cC'$ the sets of configurations which contain paths
\ts $u_1^{\sz}\lra u_{10}^{\sz}$ \ts to \ts $u_1^{\sz}\lra u_{10}^{\so}$, respectively.
The probability of a configuration \ts $Q$ and of a collection
\ts $\cC$, are given by
$$
\P(Q) \, := \, \prod_{e\ts\in\ts H \cup H'} \. p_{\psi(e)} \ts \qquad \text{and} \qquad \P(\cC) \, := \,
\sum_{Q\in \cC} \. \P(Q).
$$
We prove that \ts $\P(\cC) < \P(\cC')$ \ts by a direct combinatorial argument.

\smallskip

Note the \defna{red path \ts $\rho$} \ts from $u_1$ to $u_{10}$ in Figure~\ref{fig:h}, and
observe that it goes through every hyperedge exactly once and avoids transversal vertices.
Fix the 
order on the hyperedges of $H$ according to their appearance on the path~$\rho:$
$$(\triangleleft) \qquad (u_1,u_2,u_3), \ (u_3,u_6,u_9), \ (u_5,u_6,u_7), \ (u_2,u_4,u_5),  \ (u_4,u_7,u_8), \ (u_8,u_9,u_{10}).
$$

\smallskip

\nin
{\bf Construction:}
Let \ts $C\in \cC$ \ts be a configuration which contains a path \ts $u_1 \to u_{10}$.
The construction of the map is split into several cases.
%
%

\smallskip

\nin
\underline{{\em Case~$1$}}.  Choose $e=(abc)\in H$ to be the first hyperedge w.r.t.\
the order~$(\triangleleft)$,
such that
$$
(\psi(e), \psi(e')) \. \in \. \La^2 \quad \text{and} \quad (\psi(e),\psi(e')) \. \notin \. \Xi\ts,
$$
where
$$\aligned
\La \. & := \. \{abc, ab|c, ac|b, a|b|c\}\ts,  \quad \text{and} \\
\Xi \. & := \. \big\{\ts (abc,a|b|c)\ts, \. (a|b|c,abc)\ts , \.
(ab|c,ac|b)\ts, \. (ac|b,ab|c) \ts \big\}.
\endaligned
$$
As before, here \ts $a\in T$ \ts denotes a transversal vertex in~$e$.
To break the symmetry, we assume that $b$ precedes~$c$ along the path~$\rho$.

Note that \ts $\La^2 \ts = \ts \Om_0 \cup \Om_1 \cup \Xi$, where
$$
\Om_0 \, := \, \left\{\aligned
& \. (abc, ac|b), \. (ac|b, abc),
\. (a|b|c, ab|c), \.  (ab|c, a|b|c), \\
& \. (abc, abc), \. (a|b|c, a|b|c), \.
(ab|c, ab|c), \. (ac|b, ac|b) \. \endaligned\right\},
$$
and
\begin{equation*}
\Om_1 \, := \, \big\{\. (abc, ab|c), \.  (ab|c, abc), \.
(a|b|c, ac|b), \. (ac|b, a|b|c)\. \big\}.
\end{equation*}

\smallskip

Consider the following two cases:

\smallskip

\nin
{\small $(1)$} \. Suppose \. $(\psi(e),\psi(e'))\in \Om_0$.
Exchange the values \ts  $\psi(h) \lra \psi(h')$ \ts
for all hyperedges \emph{after} \ts $e$ or~$e'$ along the path~$\rho$.

\smallskip

\nin
{\small $(2)$} \. Suppose \. $(\psi(e),\psi(e'))\in \Om_1$.
Exchange the values \ts $\psi(e) \lra \psi(e')$, and \ts $\psi(h) \lra \psi(h')$ \ts
for all hyperedges \emph{after} \ts $e$ or~$e'$ along the path~$\rho$.

\smallskip

In both cases, denote by $C'$ the resulting configuration, and observe
that \ts $\psi(C)=\psi(C')$ \ts since we only rearrange the weights.
It remains to show the following:

\medskip

\nin
{\bf Claim:} \. $C'\in \cC'$.

\begin{proof}
Note that \ts $\psi(e) = abc$ \ts implies that vertices $b$ and~$c$ are connected
in~$C$ within~$e$ (in the WZ percolation model).
It is easy to see that in~$C'$, we have a path
$(b\to a \to a' \to c')$, i.e.\ vertices $b$ and~$c'$ are connected in~$C'$.
The same argument shows that \ts \ts $\psi(e') = abc$ \ts implies
vertices $b'$ and~$c$ are connected in~$C'$.

A similar argument shows that if vertices $a$ and~$b$ are connected in~$C$ within~$e$,
then they are still connected in~$C'$.  If vertices $a'$ and~$b'$ are
connected in~$C$ within~$e'$, then $a'$ and $b'$ are connected in~$C'$.
Similarly, if vertices $a$ and~$c$ are
connected in~$C$ within~$e$, then $a$ and $c'$ are connected in~$C'$.  Finally,
if vertices $a'$ and~$c'$ are connected in~$C$ within~$e'$, then $a$ and $c$ are
connected in~$C'$.

We have \ts $V=\{u_1,\ldots,u_{10}\}$ \ts and \ts $T=\{u_2,u_7,u_9\}$.
Let \ts $L \subseteq (V\sm T) \cup (V'\sm T')$ \ts denote the set of nontransversal
vertices which lie along the path $\rho$ between $u_1$ and~$b$, inclusively,
along with their copies from another level.  Similarly, let \ts
$R:=(V\sm T) \cup (V'\sm T') \sm L$.

Observe that for every \ts $\ga: u_1 \to u_{10}$ \ts in~$C$,
we now have a path $\ga'\in C$ constructed as follows.
For all \ts $u_i\in R$ on a path~$\ga$ use vertex $u_i'$ in~$\ga'$, and
vice versa: for all \ts $u_i'\in R$ on a path~$\ga$, use vertex $u_i'$\ts.

Recall that the values of $\psi$ are switched on all hyperedges containing vertices
in~$R$ except possibly for~$e$, and are \emph{not} switched on all hyperedges
containing vertices in~$L$ except possibly for~$e$. The result follows from
the connectivity observations above.
\end{proof}


\smallskip

\nin
\underline{{\em Case~$2$}}. In notation of Case~$1$, suppose  that \. $(\psi(e),\psi(e')) \in \Xi$.
Note that we are unable to make a switch \ts $abc \lra a|b|c$ \ts as this would disconnect
the path.  We make a \emph{probabilistic switch} defined as follows.
Denote
$$
\vk \, := \,  \. \frac{p_{ab|c} \. p_{ac|b}}{p_{abc} \. p_{a|b|c}}\.,
$$
and note that by \eqref{eq:390} we have \ts $1> \vk\ge 0$.  We make
the following switches:
$$
(ab|c,ac|b) \to (abc,a|b|c) \quad \text{and} \quad (ac|b,ab|c) \to (a|b|c,abc).
$$
Proceeding as before, we conclude that the resulting configuration
\ts $C'\in \cC'$.  However, that note we have \ts $\psi(C')-\psi(C) >0$, i.e.\
the map \ts $C\to C'$ \ts is not weight preserving.  To correct this,
for the reverse map we make switches
$$
(abc,a|b|c) \to (ab|c,ac|b)  \quad \text{and} \quad (a|b|c,abc) \to (ac|b,ab|c)
$$
with probability~$\vk$.  Proceed as in Case~1, and note that this probabilistic
map is weight preserving by construction.

\smallskip

\nin
\underline{{\em Case~$3$}}.  In notation of Case~1, we have the following
possible pairs of values of $\psi$ remain to be considered:
$$
\aligned
& (abc,a|b|c) \quad \text{and} \quad (a|b|c,abc) \quad \text{with probability \. $(1-\vk) \ts p_{abc} \. p_{a|b|c}$}\., \\
& (a|bc,\ast), \quad (\ast,a|bc)  \quad \text{and} \quad (a|bc,a|bc)  \quad \text{where} \ \ \ast \in \La\ts.
\endaligned
$$
The remaining part of the proof is a bookkeeping of this case.  The idea is that
the events in the first line emulate the alternative bunkbed hypergraph percolation
on~$\ov H$, while the events in the second line have probability too small to make
the difference by the assumption \eqref{eq:390}.

When conditioned to this case, the WZ hypergraph percolation model has the following
probabilities for each pairs of values $(\psi(e),\psi(e'))$:
$$
\begin{cases}
(abc, a|b|c), \ \text{ with probability }\frac{1}{Z} \. \big(p_{abc}p_{a|b|c} - p_{ab|c}p_{ac|b}\big),\\
(a|b|c, abc), \ \text{ with probability }\frac{1}{Z} \. \big(p_{abc}p_{a|b|c} - p_{ab|c}p_{ac|b}\big),\\
(a|bc, \ast), \ \text{ with probability }\frac{1}{Z} \. p_{a|bc} \. p_{\ast} \text{ for }\ \ast \in \La,\\
(\ast, a|bc), \ \text{ with probability }\frac{1}{Z} \. p_{a|bc} \. p_{\ast} \text{ for } \ \ast \in \La,\\
(a|bc, a|bc), \ \text{ with probability }\frac{1}{Z} \. p_{a|bc}^2\.,
\end{cases}
$$
where
$$
Z \. := \. 2\ts p_{abc} \. p_{a|b|c} \. - \.  2\ts p_{ab|c} \. p_{ac|b} \. + \.  2\ts p_{a|bc} \. - \. p_{a|bc}^2$$
is the normalizing constant.  When referring to this probabilities, we will use $\P^{(3)}$.

Denote by $\cA$ the event that for a configuration~$C$ we have \ts $\psi(e)\ne\psi(e')$ \ts
and \ts $\psi(e),\psi(e')\in \{abc, a|b|c\}$ \ts for all \ts $e\in H$.
Using the inequality \ts $(1-x)^a\ge 1-a\ts x$ \ts and the assumption \eqref{eq:390} in the lemma,
we have:
{\small
$$\aligned
\P^{(3)}(\cA) \, & =  \bigg(1 \. - \. \frac{2\ts p_{a|bc} - p_{a|bc}^2}{2(p_{abc} \.
p_{a|b|c} \ts - \ts p_{ab|c} \. p_{ac|b})  + 2\ts p_{a|bc} - p_{a|bc}^2}\bigg)^6  \ge   \left(1 - \frac{p_{a|bc}}{(p_{abc} \. p_{a|b|c} \ts - \ts  p_{ab|c} \. p_{ac|b})  +  p_{a|bc}}\right)^6\\
& \ge \, 1 \, - \, \frac{6 \ts p_{a|bc}}{(p_{abc}\ts p_{a|b|c} \. - \.  p_{ab|c}\. p_{ac|b}) \. + \. p_{a|bc}}
\, >_{\eqref{eq:390}} \, 1 \, - \, \frac{6 \ts p_{a|bc}}{401 \ts p_{a|bc}} \, > \, \frac{64}{65}\..
\endaligned
$$
}

As we mentioned above, if we condition on $\cA$, the WZ model turns into  alternative bunkbed
hypergraph percolation model, so by Hollom's result we have:
$$\aligned
\P^{(3)}(u_1^{\sz}\lra u_{10}^{\sz}\mid \cA) \. - \. \P^{(3)}(u_1^{\sz}\lra u_{10}^{\so}\mid\cA) \, & = \, \Palt(u_1^{\sz}\lra u_{10}^{\sz}) \. - \.  \Palt(u_1^{\sz}\lra u_{10}^{\sz}) \\
& =_{\text{Lemma~\ref{lem:hollom}}}
\, \tfrac{12}{64} \. - \. \tfrac{13}{64} \, = \, - \tfrac{1}{64}\..
\endaligned
$$
This implies:
\begin{equation*}
\aligned
\P^{(3)}(u_1^{\sz}\lra u_{10}^{\sz}) \. - \. \P^{(3)}(u_1^{\sz}\lra u_{10}^{\so}) \ & \le \  \P^{(3)}(\ov{\cA}) \. + \. \P^{(3)}(\cA)\Big(\P^{(3)}(u_1^{\sz}\lra u_{10}^{\sz}\mid \cA) - \P^{(3)}(u_1^{\sz}\lra u_{10}^{\so}\mid \cA)\Big) \\
&  < \  \tfrac{1}{65} - \tfrac{1}{64}\cdot\tfrac{64}{65} \ = \  0.
\endaligned
\end{equation*}
This finishes the analysis of Case~3.

\smallskip

\nin
\underline{{\em Adding the probabilities}}.
From above, for the cases 1 and~2, respectively, we have:
$$\P^{(1)}(u_1^{\sz}\lra u_{10}^{\sz}) \. - \. \P^{(1)}(u_1^{\sz}\lra u_{10}^{\so}) \ = \
\P^{(2)}(u_1^{\sz}\lra u_{10}^{\sz}) \. - \. \P^{(2)}(u_1^{\sz}\lra u_{10}^{\so}) \ = \ 0.
$$
Therefore, we have \ $\PWZ(u_1^{\sz}\lra u_{10}^{\sz}) - \PWZ(u_1^{\sz}\lra u_{10}^{\so}) < 0$,  as desired. \qed

\medskip


\section{Disproof of the Bunkbed Conjecture}\label{s:diproof}

\subsection{Hyperedge simulation}\label{ss:diproof-edge}
In this section, we construct a graph that simulates a hyperedge in the sense
of WZ hypergraph percolation, adhering to the conditions of the Lemma~\ref{lem:robust}.
We prove the following technical result for the \emph{weighted percolation}.

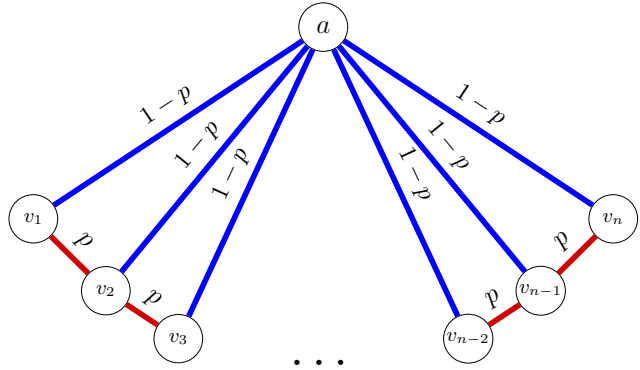
\begin{figure}[hbt]
    \centering
    \begin{tikzpicture}[scale=1.5]

\definecolor{myblue}{rgb}{0.0, 0.0, 1.0}  
\definecolor{myred}{rgb}{0.8, 0.0, 0.0} 

\tikzset{vertex/.style = {draw, circle, fill=white, text=black, minimum size=21pt, text width=21pt, align=center, inner sep=0pt, font=\footnotesize}}

\node[vertex, font=\large] (a) at (0,2) {\(a\)};
\node[vertex] (u1) at (-3,0) {$v_1$};
\node[vertex] (u2) at (-2.25,-0.75) {$v_2$};
\node[vertex] (u3) at (-1.5,-1.25) {$v_3$};
\node[vertex] (un2) at (1.5,-1.25) {$v_{n-2}$};
\node[vertex] (un1) at (2.25,-0.75) {$v_{n-1}$};
\node[vertex] (un) at (3,0) {$v_n$};

\node[scale=2.0] at (0,-1.5) {$\dots$};

\draw[myblue, line width=2.5pt] (a) -- (u1) node[midway, above, sloped, black, font=\normalsize] {$1-p$};
\draw[myblue, line width=2.5pt] (a) -- (un) node[midway, above, sloped, black, font=\normalsize] {$1-p$};
\draw[myblue, line width=2.5pt] (a) -- (u2) node[midway, above, sloped, black, font=\normalsize] {$1-p$};
\draw[myblue, line width=2.5pt] (a) -- (u3) node[midway, above, sloped, black, font=\normalsize] {$1-p$};
\draw[myblue, line width=2.5pt] (a) -- (un2) node[midway, above, sloped, black, font=\normalsize] {$1-p$};
\draw[myblue, line width=2.5pt] (a) -- (un1) node[midway, above, sloped, black, font=\normalsize] {$1-p$};

\draw[myred, line width=2.5pt] (u1) -- (u2) node[midway, above, sloped, black, font=\large] {$p$};
\draw[myred, line width=2.5pt] (u2) -- (u3) node[midway, above, sloped, black, font=\large] {$p$};
\draw[myred, line width=2.5pt] (un2) -- (un1) node[midway, above, sloped, black, font=\large] {$p$};
\draw[myred, line width=2.5pt] (un1) -- (un) node[midway, above, sloped, black, font=\large] {$p$};

\end{tikzpicture}
    \vskip-.1cm
    \caption{Graph \ts $G_n$ \ts with \ts $n+1$ \ts vertices.}
    \label{fig:gadget}
\end{figure}

\begin{lemma}\label{l:hyperedge}
Let \ts $n\ge 3$ \ts and \ts $0<p<1$.  Consider a weighted graph \ts  $G_n$ \ts on $(n+1)$ vertices
given in Figure~\ref{fig:gadget}.  Denote \ts $b:=v_1$ \ts and \ts $c:=v_{n}$.
Then \. $p_{ab|c} \. = \. p_{ac|b}$ \. and
\begin{equation}\label{eq:hyper-sim}
    p_{abc}\. p_{a|b|c} \. - \. p_{ab|c}\. p_{ac|b} \, > \, \big(n \. \tfrac{1 - p}{1 + p} \. - \. 1\big) \ts p_{a|bc}\.,
\end{equation}
where
$$\aligned
& p_{abc} \. := \. \P_p[a\lra b \lra c], \quad p_{a|bc} \. := \. \P_p[a \nleftrightarrow  b \lra c],
\quad p_{ab|c} \. := \. \P_p[a \lra  b \nleftrightarrow c], \\
& \qquad\qquad p_{ac|b} \. := \. \P_p[a \lra  c \nleftrightarrow b] \quad \text{and} \quad
p_{a|b|c} \. := \. \P_p[a \nleftrightarrow  b \nleftrightarrow c \nleftrightarrow a]\ts.
\endaligned
$$
\end{lemma}

We prove the lemma in the next section, see Proposition~\ref{prop:Gn_ineqs}.

\subsection{Proof of Theorem~\ref{thm:main}}\label{ss:diproof-proof}
In notation of Lemma~\ref{lem:robust}, let \ts $p=\frac12$ \ts and let \ts $n:= 3\cdot 401 +1 = 1204$.
The resulting graph $G_n$ is planar, has $1205$ vertices and $2407$ edges.

Take Hollom's hypergraph $H$ from Figure~\ref{fig:h} and substitute for each $3$-hyperedge
with a graph \ts $G_n$ \ts from Lemma~\ref{l:hyperedge}, placing it so $a$ is a transversal vertex
while \ts $b=v_1$ \ts and \ts $c=v_n$ \ts are the other two vertices.
The resulting graph is still planar, has \.
$10 + 6\cdot 1202 = 7222$ \. vertices and \. $6 \cdot 2407 =14442$ \. edges.

By Lemma~\ref{l:hyperedge}, the \ts $\frac12$-percolation on
\ts $G_n$ \ts satisfies conditions of Lemma~\ref{lem:robust}.   Thus, by
Lemma~\ref{lem:robust}, we have:
$$\P(u_1^{\sz}\lra u_{10}^{\sz}) \ < \ \P(u_1^{\sz}\lra u_{10}^{\so}),
$$
as desired.
\qed

\begin{rem}\label{r:counter-example}
Due to the multiple conditionings and the gadget structure, the difference of probabilities
given by the counterexample in Theorem~\ref{thm:main}
is less than \ts $10^{-4331}$, out of reach computationally.
%
%
A computer-assisted computation shows that one can use \ts $G_n$ \ts with
\ts $p = \frac12$ \ts  and \ts $n=14$, giving a relatively small graph on $82$ vertices.
However, even in this case, the difference of the probabilities in the BBC
is on the order \ts $10^{-47}$.  This and other computations are collected
on author's website, see~$\S$\ref{ss:finrem-lemma}.

Since \emph{Weighted BBC} is equivalent to BBC (see~$\S$\ref{ss:finrem-models}),
once can instead take weighted graph \ts $G_n$ \ts with \ts $p = \frac{1}{2n}$ \ts and \ts $n=402$.
This graph is still too large to analyze experimentally.
A computer-assisted computation shows that one can use \ts $G_n$ with \ts $p = 0.0349$ \ts
and \ts $n=5$, giving a rather small graph on $28$ vertices.  However, even in this case,
the difference of the probabilities in the Weighted BBC are on the order \ts $10^{-78}$.
\end{rem}

\medskip

\section{Proof of Lemma~\ref{l:hyperedge}}
\label{s:Gn}

We prove the lemma as a consequence of elementary calculations.

\smallskip

\begin{lemma}\label{lem:p_ab}
We have:
$$\P_p(a \leftrightarrow v_n) \, =  \, \frac{1 - p^{2n}}{1 + p}\..$$
 \end{lemma}

\begin{proof}
Let \ts $p_n :=\P_p(a \leftrightarrow v_n)$ \ts as in the lemma.
We establish a recurrence relation for~$p_n$. There are two cases:

\smallskip

\nin
{\small $(1)$} \. The edge \ts $(a, v_n)$ is open. This occurs with probability \ts $1 - p$. In this case,
vertices \ts $a$ \ts and \ts $v_n$ \ts are directly connected.

\smallskip

\nin
{\small $(2)$} \. The edge \ts $(a, v_n)$ \ts is closed. This occurs with probability~$p$. In this case, vertex \ts $v_n$ \ts
can only connect to \ts $a$ \ts through the edge \ts $(v_{n-1}, v_n)$, which is open with probability~$p$.
If this edge is closed, vertex \ts $v_n$ is isolated from \ts $a$.
If it is open, the probability that $a$ and $v_{n-1}$ are in the same connected component is \ts $p_{n-1}\ts$.

\smallskip

Combining these cases, we obtain the following recurrence relation:
    \[
    p_n \, = \, (1 - p) \. + \. p^2 \ts p_{n-1}\.,
    \]
    with the initial condition \ts $p_0 = 0$.  The result follows by induction.
\end{proof}

\smallskip

\begin{lemma}\label{lem:p_abc}
We have:
$$\P_p(a \leftrightarrow v_1 \leftrightarrow v_n) \ = \ \frac{1 - p^{2n}}{(1 + p)^2} \, + \,  \frac{n \ts (1-p) \ts p^{2n-1}}{1+p} \..
$$
\end{lemma}

\begin{proof}
Let \ts $p_n := \P_p(a \leftrightarrow v_1 \leftrightarrow v_n)$ \ts denote the probability as in the lemma.
We calculate this probability by analyzing whether edges \ts $(a, v_1)$ \ts and \ts $(a, v_n)$ \ts
are open or closed.  There are four cases:

\smallskip

\nin
{\small $(1)$} \. Both edges $(a, v_1)$ and $(a, v_n)$ are open, each with probability $1 - p$.
Then $a$ is directly connected to both $v_1$ and $v_n$. Thus, the probability is \ts $(1 - p)^2$.

\smallskip

\nin
{\small $(2)$} \. Edge \ts $(a, v_n)$ \ts is closed. If the edge $(a, v_n)$ is closed,
 vertex \ts $v_n$ \ts is connected to the rest of the graph through the edge \ts $(v_{n-1}, v_n)$,
 which is open with probability~$p$.  This reduce the problem to \ts $G_{n-1}$\ts.
 Thus, the probability is \ts $p^2 \ts p_{n-1}$.

\smallskip

\nin
{\small $(3)$} \. The edge $(a, v_1)$ is closed.  Similarly, if the edge \ts $(a, v_1)$ \ts
is closed (with probability $p$).  Thus, the probability is \ts $p^2 \ts p_{n-1}$.

\smallskip

\nin
{\small $(4)$} \. Both edges $(a, v_1)$ and $(a, v_n)$ are closed. If both edges $(a, v_1)$ and $(a, v_n)$ are closed (each with probability $p$), $v_1$ must connect to $v_2$ by the edge $(v_1, v_2)$, and $v_n$ must connect to $v_{n-1}$ by the edge $(v_{n-1}, v_n)$. The problem reduces to finding the probability that $a$, $\widehat{u}_1 = v_2$, and $\widehat{u}_{n-2} = v_{n-1}$ are in the same connected component in the graph $\widehat{G}_{n-2}$, the restriction of $G_n$ to the vertices $a, v_2, \dots, v_{n-1}$. Thus, the corresponding probability is \ts $p^4 \ts p_{n-2}$.

\smallskip

Using inclusion-exclusion of these four cases, we obtain the following recurrence relation:
    \[
    p_n \, = \, (1 - p)^2 \. + \. 2 \ts  p^2 \ts p_{n-1} \. - \. p^4 \ts p_{n-2},
    \]
    with initial conditions \ts $p_0 = 0$ \ts and \ts $p_1 = 1 - p$.  The result
    follows by induction.
 \end{proof}

\smallskip

\begin{lemma}\label{lem:p_a_bc}
We have:
    \[
    \P_p(a \nleftrightarrow v_1 \leftrightarrow v_n) \, = \, p^{2n - 1}.
    \]
\end{lemma}

\begin{proof}
    If the vertices $v_1$ and $v_n$ are in the same connected component that does not contain vertex~$a$,
    they must be connected by the path \. $\ga := (v_1 \to v_2 \to \. \ldots \. \to v_n)$. The probability that
    this path is open is \ts $p^{n - 1}$. In addition, any edge \ts $(a, v_k)$ \ts must be closed for all \ts
    $1 \le k \le n$, as otherwise vertex \ts $a$ \ts is connected to the path~$\ga$. The probability that
     all these edges are closed is $p^{n}$. Thus, the probability in the lemma
    is \ts $p^{2n - 1}$.
\end{proof}

\smallskip

We conclude with the following result which immediately implies Lemma~\ref{l:hyperedge}.

\smallskip

\begin{prop}\label{prop:Gn_ineqs}
In notation of Lemma~\ref{l:hyperedge}, we have \. $p_{a|bc} = p^{2n-1}$ \. and
    \begin{align*}
        p_{abc} \. p_{a|b|c} \. - \. p_{ac|b} \. p_{ab|c} \, \ge \, \left(n \. \tfrac{1-p}{1+p} \. - \. 1\right) p^{2n - 1}.
    \end{align*}
\end{prop}

\begin{proof}
The first part is given by Lemma~\ref{lem:p_a_bc}.  For the second part, using
Lemmas~\ref{lem:p_ab}, \ref{lem:p_abc} and~\ref{lem:p_a_bc} and \ts $p_{abc}\le 1$, we have:
\begin{align*}
p_{abc}\. p_{a|b|c} \. - \. p_{ac|b} \. p_{ab|c} \ &= \ p_{abc} \. - \. (p_{abc} +
p_{ab|c})(p_{abc} + p_{ac|b}) \. - \. p_{abc}\. p_{a|bc} \\
&= \ \P_p(a \leftrightarrow v_1 \leftrightarrow v_n) \. -  \.
\P_p(a \leftrightarrow v_1) \cdot \P_p(a \leftrightarrow v_n) \. - \. p_{abc}\. p_{a|bc}\\
&\ge \ \left(\tfrac{1 - p^{2n}}{(1 + p)^2} \. + \. \tfrac{n\ts (1-p) \ts p^{2n-1}}{1+p}\right) \, - \, \left(\tfrac{1 - p^{2n}}{1 + p}\right)^2  - p^{2n - 1}\\
&\ge \ \tfrac{p^{2n}(1 - p^{2n})}{(1 + p)^2} \. + \. \tfrac{n\ts (1-p) \ts p^{2n-1} }{1+p}  \. - \. p^{2n - 1}\\
& \ge \ \left(\tfrac{n\ts (1-p)}{1+p} - 1\right)p^{2n - 1},
\end{align*}
as desired.
\end{proof}

\medskip

\section{Complete BBC}
\label{s:complete}

In notation of the Bunkbed Conjecture~\ref{conj:BB}, one can ask
if a version of the BBC holds for uniform \ts $T\subseteq V$.  This
is equivalent to $\frac12$-percolation on the product graph $G\times K_2$\ts.
To distinguish from BBC, we call this \defn{Complete BBC}, see~$\S$\ref{ss:finrem-models}.
Turns out the proof of Theorem~\ref{thm:main} extends to the proof of Complete BBC,
but a counterexample is a little larger:

\smallskip

\begin{thm}\label{thm:complete}
There is a connected graph \ts $G=(V,E)$ \ts \ts with  \ts $|V|= 7523$ \ts vertices and \ts
$|E|=15654$ \ts edges, and vertices \ts $u,v \in V$, s.t.\ for the \ts
$\frac12$-percolation on \ts $G\times K_2\ts$ we have:
\[\P_{\frac12}[\ts u \lra v \ts ] \. < \. \P_{\frac12}[ \ts u \lra v' \ts ]\ts.
\]
In particular, the complete bunkbed conjecture is false.
\end{thm}

\begin{proof}[Sketch of proof]
Recall that Hollom's Model~4.3 in \cite{Hol24} is the hypergraph
version of the Complete BBC.  Hollom disproves it in \cite[$\S$5.1]{Hol24}
by showing that his $3$-hypergraph in Figure~\ref{fig:h} is minimal in a
sense that bunkbed probabilities are equal for all subsets \ts
$\{u_2,u_7,u_9\} \ssu T\subseteq \{u_1,\ldots,u_{10}\}$.  He then
makes \ts $k=102$ \ts ``clones'' of vertices \ts $\{u_2,u_7,u_9\}$ \ts to make
sure at least one is always in the percolation cluster with high
probability.

From this point on, proceed as in the proof of Theorem~\ref{thm:complete}.
In notation of the proof of Lemma~\ref{lem:robust}, we have that the
only way we can have a nonzero probability gap if the path~$\ga$ goes along
the red path~$\rho$.  In notation of the proof of Lemma~\ref{l:hyperedge},
having posts at all $v_i$ along the \ts $b\to c$ \ts paths within a
hyperedge \ts $(abc)$ \ts gives zero probability gap unless all posts
are closed.  Combining these observations proves the theorem.
The resulting graph has \ts $|V|=7222 + 3(k-1)= 7523$ \ts and \ts
$|E|=14442 + 3\cdot 4\cdot (k-1) = 15654$.  We omit
the details.  \end{proof}

\begin{rem}\label{r:complete}
In notation of Remark~\ref{r:counter-example}, the difference of
probabilities is even smaller in this case, and is on the
order of $10^{-6509}$.  This is a consequence of having all
posts that are not clones of transversal posts being closed.
\end{rem}

\medskip

\section{Experimental testing}
\label{s:experimental}

Versions of the bunkbed conjecture were repeatedly tested by various
researchers, although failed attempts remain largely unreported,
see e.g.\ \cite[$\S$3.1]{Rud21}.  In this section we describe our
own attempt to refute the conjecture using a large scale computer
computation.

\subsection{Initial tests}\label{ss:exp-first}
We started with a series of brute force testing
of the Polynomial BBC, see~$\S$\ref{ss:finrem-models}.
We exhaustively tested all connected graphs with at most
8~vertices, and connected graph with at most 15 edges
from the \ts {\tt House of Graphs}  \ts database, see \cite{CDG23}.
In each case, the Polynomial BBC held true.  At this point
we chose to develop a more refined approach.

\subsection{The algorithm}\label{ss:exp-algorithm}
Our starting point is the machine learning algorithm by
Wagner \cite{Wag21}, which we adjusted and modified.
Roughly, the algorithm inputs a neural network used in a randomized graph
generating algorithm, various constraints and a function to optimize.
It outputs new weights for the neural network with the function improved.
In his remarkable paper, Wagner describes how he was able to disprove
five open problems in graph theory, so we had high hopes that this
approach might help to disprove the bunkbed conjecture.

To give a quick outline of our approach, we consider a graph \ts $G=(V,E)$
on $n=|V|$ vertices, with the set of transversal vertices \ts $T\ssu V$,
and fixed \ts $u,v\notin W$. Flip a fair coin for each edge $e\in E$.
Depending on the outcome, either retain~$e$ and delete~$e'$,
or vice versa.  Check whether \ts $u \lra v$ \ts and \ts $u \lra v'$.
Repeat this $N$ times to estimate the corresponding probabilities
$p$ and~$p'$, respectively.  Based
on these probabilities, use Wagner's algorithm to obtain the
next iteration.  Repeat $M$ times or until a potential
counterexample with $p< p'$ is found.

\subsection{Implementation and results}\label{ss:exp-implementation}
We first used Wagner's original code on a laptop computer,
but when that proved too slow we made major changes.  To speed
up the performance and tweak the code, we implemented Wagner's
algorithm in {\sc Julia}.

We then run the code on a shared {\sc UCLA Hoffman2 Cluster}, which is a
Linux compute cluster consisting of 800+ 64-bit nodes and over 26,000 cores,
with an aggregate of over 174TB of memory.\footnote{The system overview
is available here: \ts
\url{www.hoffman2.idre.ucla.edu/About/System-overview.html}}
Each run of the algorithm required about 2 hours.  After six runs
with different parameters, the results were too similar to continue.

Specifically, we run the algorithm on graphs with $n=20$ and $n= 30$ vertices,
and for~3, 4 and 5 transversal vertices.  Although we started with relatively
dense graphs, the algorithm converged to relatively sparse graphs with about
100 edges. We used $N=4000$, pruning the Monte Carlo sampling when the
desired probabilities were far apart.

We used $M=2000$, after which the probabilities $p,p'$ rapidly converged to $\frac12$
and became nearly indistinguishable.   More precisely, the \emph{probability gap}
\ts $p-p'$ \ts became smaller than $0.01$ getting close to the Monte Carlo error,
i.e.\ the point when we would need to increase $N$ to avoid false positives.
At all stages, we had \ts $p> p'$ \ts suggesting validity of the bunkbed conjecture.
At the time of the experiments and prior to this work, we saw no evidence
that an experimental approach could ever succeed.

\subsection{Analysis}\label{ss:exp-ana}
Having formally disproved the bunkbed conjecture,
it is clear that our computational approach was doomed to failure.
For the uniform weights we tested, we could never have reached graphs
of size anywhere close to that in Theorem~\ref{thm:main}, of course.
Even when the number of vertices is optimized to 82 as suggested in
Remark~\ref{r:counter-example}, the number of edges is still
very large while the probability gap in the theorem is
on the order of \ts $10^{-47}$, thus undetectable in practice.

In hindsight, to reach a small counterexample we should have used
the weighted bunkbed percolation rather than the more efficient
alternative model, with some edges having a very large weight and
some very small weight.  Of course, by Remark~\ref{r:counter-example},
the probability gap in the theorem is still prohibitively small,
at least for the graphs we consider.

\subsection{Conclusions}\label{ss:exp-con}
It seems, the Bunkbed Conjecture has some unique features
making it very poorly suited for computer testing. In fact,
one reason we stopped our computer experiments is that in our initial
rush to testing we failed to contemplate societal implications
of working with even moderately large graphs.

Suppose we \emph{did} \ts find a
potential counterexample graph with only \ts $m=100$ \ts edges and the
probability gap was large enough to be statistically detectable.
Since analyzing all of \ts $2^m \approx 10^{30}$ \ts subgraphs is
not feasible, our Monte Carlo simulations could only confirm
the desired inequality with high probability.  While this
probability could be amplified by repeated testing, one could
never formally \emph{disprove} the bunkbed conjecture this way,
of course.

This raises somewhat uncomfortable questions whether the
mathematical community is ready to live with an uncertainty over
validity of formal claims that are only known with high probability.
It is also unclear whether in this imaginary world the granting
agencies would be willing to support costly computational projects
to further increase such probabilities (cf.~\cite{G+,Zei93}).
Fortunately, our failed computational effort avoided this
dystopian reality, and we were able to disprove the bunkbed
conjecture by a formal argument.

\medskip

\section{Final remarks}\label{s:finrem}

\subsection{Variations on the BBC} \label{ss:finrem-models}
Although the version of the Bunkbed Conjecture~\ref{conj:BB} given in \cite{BKa}
is considered the most definitive, over the years several closely related
versions has been studied:
\smallskip

\nin
{\small $(0)$} \. \defn{Counting BBC}, where one compares the number of subgraphs
connecting vertices $u,v$ and those that do not.  This conjecture is a restatement
of the BBC in the case \ts $p=\frac12$.

\smallskip

\nin
{\small $(1)$} \. \defn{Weighted BBC}, where the edge probabilities \ts $p_e=p_{e'}$ \ts
can depend on $e\in E$.  This conjecture is equivalent to the BBC by \cite{RS},
since edge probabilities can be approximated by series-parallel graphs.

\smallskip

\nin
{\small $(2)$} \. \defn{Polynomial BBC}, where the edge probabilities above are viewed
as variables.  In this case the conjecture claims that the difference of polynomials
on is a polynomial in nonnegative coefficients.  This conjecture is stronger than
Weighted BBC as there are polynomials positive on $[0,1]$ which have negative
coefficients such as \ts $(x-y)^2$.  Although we did not find a counterexample
on graph with at most 8 vertices, is likely that there is a sufficiently small
counterexample in this case.  Cf.\ \cite{Ric22}, where the difference is a sum
of squares.

\smallskip

\nin
{\small $(3)$} \. \defn{Computational BBC}, where one asks if the counting version
of the probability gap is in $\SP$, i.e.\ has a combinatorial interpretation,
see \cite[Conj.~5.6]{Pak-OPAC}.  Clearly, this conjecture
implies BBC.  We refer to \cite{IP22} for a formal treatment
of this problem for general polynomials.

\smallskip

\nin
{\small $(4)$} \. \defn{Alternative BBC}, where fair coin flips determine whether
the edge $e$ is deleted and $e'$ retains or vice versa.  This conjecture
implies BBC \cite[Prop.~2.6]{Lin11}.

\smallskip

\nin
{\small $(5)$} \. \defn{Complete BBC}, where one takes all $T=V$ and performs
the weighted percolation on the full \ts $\ov G := G\times K_2$\ts, i.e.\ on
all edges in $\ov G$ including the posts.  The conjecture in this case is weaker
than the BBC, see e.g.\ \cite[Prop.~2.2]{Lin11}.

\smallskip

In all but the last case, the corresponding conjecture is refuted by Theorem~\ref{thm:main}.
In {\small $(5)$}, the corresponding conjecture is refuted by Theorem~\ref{thm:complete}
using the same counterexample.

\subsection{Robustness lemma}\label{ss:finrem-lemma}
Lemma~\ref{lem:robust} is a finite problem which can be reformulated as follows.
By definition, probabilities on both sides of \eqref{eq:semi} are polynomials
in 5 variables of degree at most~12, with at most $5^{12}$ nonzero coefficients.
The lemma gives positivity of the difference of these two polynomials on a
region of the unit cube \ts $[0,1]^5$ \ts cut out by the quadratic
inequality~\eqref{eq:390}.

Since our proof of Lemma~\ref{lem:robust} is somewhat cumbersome and uses
a case-by-case analysis, we verified the lemma computationally.  The results
and the code are available on {\tt GitHub}.\footnote{\emph{Generating Partitions of Graph Vertices into Connected Components},
description and  
 code at \\ \href{https://github.com/Kroneckera/bunkbed-counterexample}{GitHub.com/Kroneckera/bunkbed-counterexample}}
Of course, the advantage of our combinatorial proof is
that it is elementary and amenable for generalizations.

\subsection{Special cases}\label{ss:finrem-special}
Our counterexample makes prior positive results somewhat more valuable.
It would be interesting to find other families of graphs on which the
Bunkbed Conjecture holds.  We are especially interested in the
corresponding problem for the Polynomial BBC.  Note that we
emphasized planarity in Theorem~\ref{thm:main} since it was
speculated in \cite{Lin11} that planarity helps.

\subsection{Two-point functions}\label{ss:finrem-2PF}
For lattices, the connection probabilities \ts $\P_p[u\lra v]$ \ts
between vertices are known as the \emph{two-point functions}.
For percolation in higher dimensions, these were famously studied
by Hara, van der Hofstad and Slade \cite{HHS03}, and they are
also of interest for other probabilistic models.
Curiously, it is not known whether connection probabilities
are monotone as the distance \ts $|u-v|$ \ts increases.
This claim would follow from the bunkbed
conjecture.  This suggests that investigating the BBC
for grid-like graphs is still of interest even if the
conjecture is false for general planar graphs.
Note that the monotonicity is known in the \ts $p\downarrow 0$ \ts limit,
and for the Ising model \cite{MMS77}.

\subsection{Random cluster model}\label{ss:finrem-RC}
It was shown in \cite[$\S$3]{Hag03} that the analogue of the BBC holds
for the random cluster model with parameter \ts $q=2$.  Our
Theorem~\ref{thm:main} shows that one cannot take \ts $q=1$.
It would be interesting to find the smallest \ts $q>1$ \ts
s.t.\ the BBC holds for all finite graphs.

\vskip.6cm
{\small
\subsection*{Acknowledgements}
Soon after joining UCLA in 2009, IP learned about the bunkbed conjecture from Tom Liggett,
who strongly encourages us to work on the problem.  We are saddened that Tom did not get
to see the disproof of the conjecture, and dedicate the paper to his memory.

We are grateful to B\'{e}la~Bollob\'as, Swee Hong Chan, Tom Hutchcroft, Jeff Kahn,
Piet Lammers, Svante Linusson, Bhargav Narayanan, Petar~Nizi\'c-Nikolac and
Fedya Petrov for helpful discussions and remarks on the subject, and to
Doron Zeilberger for jovial comments.  We thank Lisa Snyder from the
UCLA Advanced Research Computing team, for helping us with the
{\sc Hoffman2 Cluster} setup.

This paper was written when IP was a member at the Institute of Advanced Study (IAS)
in Princeton, NJ, and was partially supported by the grant from the IAS School of Mathematics;
we are grateful for the hospitality.  NG and IP were partially supported by
the NSF CCF Award No.~2302173.
}

\vskip.9cm


\end{document}